\documentstyle[12pt,amssym]{article}
\newtheorem{theo}{Theorem}

\newtheorem{lemm}{Lemma}
\newtheorem{coro}{Corollary}

\newcommand{\myF}{{\cal F}}

\newcommand{\mya}{{\frak a}}
\newcommand{\myg}{{\frak g}}
\newcommand{\myh}{{\frak h}}
\newcommand{\myk}{{\frak k}}

\newcommand{\mym}{{\frak m}}

\newcommand{\myp}{{\frak p}}

\newcommand{\bR}{{\Bbb R}}

\newcommand{\bX}{{\Bbb X}}
\newcommand{\bY}{{\Bbb Y}}
\newcommand{\Ad}{{\rm Ad}}

\newcommand{\trace}{{\rm Trace}}
\newcommand{\diag}{{\rm diag}}
\newcommand{\id}{{\rm id}}

\newcommand{\levy}{L\'{e}vy }

\begin{document}

\begin{center}
{\large \bf A decomposition of Markov processes via group actions}

Ming Liao\footnote{Department of Mathematics, Auburn University,
Auburn, AL 36849, USA. Email: liaomin@auburn.edu}
\end{center}

\begin{quote}
{\bf Summary} \ We study a decomposition of a general Markov process in a
manifold invariant under a Lie group action into a radial part (transversal to orbits)
and an angular part (along an orbit). We show that given a radial path,
the conditioned angular part is a nonhomogeneous \levy process in
a homogeneous space, we obtain a representation of such processes, and
as a consequence, we extend the well known skew-product of Euclidean Brownian motion
to a general setting.


\end{quote}

\section{Introduction} \label{intro}

It is well known that a Brownian motion $x_t$ in $\bR^n$ ($n\geq 2$) may be expressed
as a skew product of a Bessel process and a spherical
Brownian motion. The Bessel process $r_t=\vert x_t\vert$ is the radial part of $x_t$
and the angular part $\theta_t=x_t/r_t$
is a timed changed spherical Brownian motion.
This decomposition is naturally related to the action of the rotation group $SO(n)$,
the group of $n\times n$ orthogonal matrices of determinant $1$,
as the Brownian motion $x_t$ has an $SO(n)$-invariant distribution with its radial part
transversal to the orbits of $SO(n)$ and angular part contained in an orbit,
namely the unit sphere. More generally, it is shown in Galmarino \cite{galmarino} that
a continuous Markov process in $\bR^n$ with an $SO(n)$-invariant distribution
is a skew product of its radial motion and an independent spherical Brownian motion with
a time change. Such a skew-product structure in connection with a group action has been
noticed in literature. For example, Pauwels and Rogers \cite{PR} considered
a skew-product of Brownian motion in a manifold.

In this paper, we will consider a general Markov process $x_t$ in a smooth manifold $X$
that has a distribution invariant under the smooth action of a Lie group $K$.
Let $Y$ be a submanifold of $X$ transversal to the orbits of $K$. The radial part and
the angular part of $x_t$ are respectively its projections to $Y$ and to a typical $K$-orbit.
It is easy to show that the radial part is a Markov process in $Y$. Our main purpose is
to study the conditioned angular process given a radial path.

In the next section, we provide the easy proof that the radial part $y_t$ of $x_t$
is a Markov process in $Y$. We also briefly discuss to what extend the process $x_t$
is determined by its radial part. Several examples are mentioned here.

In order to study the angular part, we need first to introduce, in section~3,
the notion of nonhomogeneous \levy processes (processes of independent
but not necessarily stationary increments) in a homogeneous space, and
establish their representation in terms of a drift, a covariance operator and
a \levy measure function. This is an extension of Feinsilver's result \cite{feinsilver} for
nonhomogeneous \levy processes in Lie groups. The arguments
in \cite{feinsilver} may be suitably modified to work on homogeneous spaces, but requires
a careful formulation of a product structure on homogeneous spaces.

The angular part of $x_t$ is introduced in section~4 as a process in a typical $K$-orbit
which may be identified with a homogeneous space $K/M$.
We prove that given a radial path, the conditioned angular process $z_t$
is a nonhomogeneous \levy process in $K/M$. As a consequence, we show that
if $x_t$ is continuous and if $K/M$ is irreducible, then $x_t$ is a skew product of
its radial part and an independent Brownian motion in $K/M$ with a time change.
This is an extension of Galmarino's result to a more general setting
by a conceptually more transparent proof.

In section~5, we study a class of $K$-invariant Markov processes in $X$ obtained
by interlacing a diffusion process with jumps. In this case, we may obtain explicit
expressions for the covariance operator and the \levy measure function of
the conditioned angular process.

The skew-product decomposition of Brownian motion in Euclidean spaces or manifolds
may be well known, but this is perhaps the first time when it is studied in the general setting
of a Markov process in a manifold under the action of a Lie group, assuming only
a simple Markov property with a possibly finite life time and some technical condition
on the group action.
We like to mention that it may be possible to extend our results to the action of
a locally compact group $K$, using an extension of the representation of nonhomogeneous
\levy processes in such groups, obtained in Heyer and Pap \cite{hp}, to
homogeneous spaces.

\section{Radial part of a Markov process} \label{radial}

Throughout this paper, let $X$ be a (smooth) manifold and let $K$ be a Lie group
acting (smoothly) on $X$. Let $C_b(X)$, $C_c(X)$ and $C_0(X)$ be the spaces of continuous
functions on $X$ that are respectively bounded, compactly supported and convergent
to $0$ at infinity in the one-point compactification of $X$. When a superscript $\infty$
is added, such as $C_c^\infty(X)$, it will denote the subspace of smooth functions.

Let $x_t$ be a Markov process in $X$ with transition semigroup $P_t$.
By this we mean a process $x_t$ with rcll paths
(right continuous paths with left limits) that has the following simple
Markov property:
\begin{equation}
E[f(x_{t+s})\mid\myF_t] = P_sf(x_t) \label{simpleMP}
\end{equation}
almost surely for $s<t$ and $f\in C_b(X)$,
where $\myF_t$ is the natural filtration of process $x_t$, and for $t\geq 0$, $P_t$
is a sub-probability kernel from $X$ to $X$ (that is, $P_t(x,\cdot)$ is a sub-probability
measure on $X$ for $x\in X$ and $P_t(x,B)$ is measurable in $x$ for measurable
$B\subset X$), with $P_0(x,\cdot)=\delta_x$ (unit point mass at $x$), such that
$P_{t+s}(x,\cdot)=\int P_t(x,dy)P_s(y,\cdot)$. Note that the Markov process $x_t$
is allowed to have a finite life time as $P_t(x,X)$ may be less than $1$.

We will assume the Markov process $x_t$ or equivalently its transition semigroup $P_t$
is $K$-invariant in the sense that
\begin{equation}
\forall f\in C_b(X)\ {\rm and}\ k\in K,\ \ \ \ P_t(f\circ k) = (P_tf)\circ k.
\label{Kinv}
\end{equation}
This means that for $k\in K$, $kx_t$ is the same Markov process started at $kx_0$
(in the sense of distribution).

The Markov process $x_t$ in $X$ is called a Feller process if $P_tf\in C_0(X)$
for $f\in C_0(X)$ and $P_tf\rightarrow f$ uniformly as $t\rightarrow 0$.
In this case, $P_t$ is completely determined by its generator $L$ given by
$Lf=\lim_{t\to 0}(1/t)P_tf$ with domain $D(L)$ consisting of $f\in C_0(X)$
for which the limit exists under the sup norm.
A continuous Feller process in $X$ will be called a diffusion
process in $X$ if its generator $L$ restricted to $C_c^\infty(X)$
is a differential operator with smooth coefficients that annihilates constants.

Let $Y$ be a submanifold of $X$, possibly with a boundary, that is transversal
to the action of $K$ in the sense that it intersects each orbit of $K$ at exactly
one point, that is,
\begin{equation}
\forall y\in Y,\ \ \ \ (Ky)\cap Y=\{y\}\ \ \ \ {\rm and}
\ \ \ \ X = \cup_{y\in Y}Ky. \label{transv}
\end{equation}
Let $J$: $X\rightarrow Y$ be the projection map $J(x)=y$ for $x\in Ky$, which
is continuous if $K$ is compact. Note that $J\circ k=J$ for $k\in K$.

\begin{theo} \label{th1}
$y_t=J(x_t)$ is a Markov process in $Y$ with transition semigroup $Q_t$ given by
\begin{equation}
Q_t f(y)=P_t(f\circ J)(y),\ \ \ \ y\in Y\ \ {\rm and}\ \ f\in C_b(Y).
\end{equation}
Moreover, if $x_t$ is a Feller process in $X$ with generator $L$ and if $K$ is compact,
then so is $y_t$ in $Y$ with generator $L^Y$ given by $(L^Yf)\circ J=L(f\circ J)$
for $f\in D(L^Y)=\{h\circ J$; $h\in D(L)\}$.
\end{theo}

\noindent {\bf Proof} \ For $f\in C_b(Y)$ and $y\in Y$,
\begin{eqnarray*}
&& E[f(y_{t+s})\mid\myF_t] = E[f\circ J(x_{t+s})\mid\myF_t] = P_s(f\circ J)(x_t) \\
&=& P_s(f\circ J\circ k^{-1})(kx_t)\ \ \ \ \mbox{(where $k\in K$ is chosen such that
$kx_t=y_t$)} \\
&=& P_s(f\circ J)(y_t).
\end{eqnarray*}
This proves that $y_t$ is a Markov process in $Y$ with transition semigroup $Q_t$.
If $K$ is compact, $f\circ J\in C_0(X)$ for $f\in C_0(Y)$, the Feller property of $y_t$
follows from that of $x_t$. \ $\Box$
\vspace{2ex}

The process $y_t=J(x_t)$ in Theorem~\ref{th1} will be called the radial part of
process $x_t$ (relative to $K$ and $Y$). Note that for a diffusion process $x_t$
with generator $L$, the generator $L^Y$ of $y_t$ is the radial part of the differential
operator $L$ as defined in \cite{helgason2}.

For a measure $\mu$ and a function $f$ on $X$, the integral $\int f(x)\mu(dx)$ maybe
written as $\mu(f)$. For a measurable map $g$: $X\to X$, let $g\mu$ be the measure
on $X$ defined by $g\mu(f)=\mu(f\circ g)$.

If $K$ is compact and the distribution $\mu$ of $x_0$ is $K$-invariant,
that is, if $k\mu=\mu$ for $k\in K$, then the marginal distributions of $x_t$
are completely determined by those of the radial part $y_t$.
Indeed, For $f\in C_b(X)$ and $k\in K$,
\[E[f(x_t)] = \mu(P_tf) = \mu[(P_tf)\circ k] = \mu[P_t(f\circ k)] =
\mu(P_t\tilde{f}) = J\mu(Q_t\tilde{f}) = E[\tilde{f}(y_t)],\]
where $\tilde{f}=\int dk(f\circ k)$, with $dk$ being the normalized Haar measure on $K$,
may be regarded either as a $K$-invariant function on $X$ or a function on $Y$.

However, the distribution of the process $x_t$ is not in general determined by
the distribution of its radial part $y_t$. Consider a process in $\bR^n$, when starting
at a point different from the origin, it is a Bessel process in the ray containing
the starting point, but
when starting at the origin, it immediately chooses a ray with a uniform distribution
and then performs Bessel motion along the ray. It is easy to see that this is
a continuous Feller process in $\bR^n$ that has the same radial part as a Brownian motion
in $\bR^n$, but not a Brownian motion in distribution.

The radial part determines the process $x_t$ if the transition semigroup $P_t$
is determined by $Q_t$, which will require an additional assumption.
For example, let $G$ be a Lie group acting transitively on $X$ and containing $K$
as a subgroup, and assume there is a point $o$ in $Y$ fixed by $K$.
Then the transition semigroup $P_t$ of a $G$-invariant Markov process in $X$
is determined by $Q_t(o,\cdot)$, where $Q_t$ is the transition semigroup of the radial part.
To prove this, let $x_t$ start at $o$. Then $P_t(o,\cdot)$ is determined
by $Q_t(o,\cdot)$. Now the $G$-invariance of $P_t$ and the transitivity
of the $G$-action on $X$ imply that $P_t$ is determined by $P_t(o,\cdot)$
and hence by $Q_t$.

A Feller process $x_t$ in a Lie group $G$ with an infinite life time,
invariant under left translations, will be called
a \levy process in $G$. Such a process possesses independent and stationary increments in
the sense that for any $s<t$, $x_s^{-1}x_t$ is independent of the process up to time $s$
and its distribution depends only on $t-s$ (see \cite{liao}). This notion extends
the usual definition of \levy processes in $\bR^n$ regarded as
an additive group.

More generally, a Feller process $x_t$ in a homogeneous space $G/K$,
with an infinite life time and invariant under the natural (left) action of $G$ on $G/K$, 
where $K$ is a compact subgroup of $G$, will also
be called a \levy process in $G/K$. It is clearly $K$-invariant,
and by the discussion in the last paragraph, its distribution is determined by
its radial part relative to any $Y$ transversal to $K$ that has a point fixed by $K$.

An explicit formula for the generator of a \levy process in $G$ or $G/H$ is obtained
by Hunt \cite{hunt}, see also \cite[chapter 2]{liao}.
A \levy process in $G$ may be characterized by a stochastic integral
equation driven by a Brownian motion and a Poisson random measure, see Applebaum
and Kunita \cite{ak}.
\vspace{2ex}

For $x\in X$, $K_x=\{k\in K$; $kx=x\}$ is a closed subgroup of $K$, called
the isotropy subgroup of $K$ at $x$. Let $Y^\circ$ be $Y$ minus its boundary.
Let $X^\circ$ be the union of the $K$-orbits that intersect $Y^\circ$. This is an
open dense subset of $X$. In order to introduce an angular part of
the Markov process $x_t$ later in section~\ref{angular}, we will assume that $K_y$
is the same compact subgroup $M$ of $K$ as $y$ varies over $Y^\circ$.
This assumption is often satisfied when the transversal submanifold $Y$ is properly chosen.
Then $X^\circ=Y^\circ\times(K/M)$ as a product manifold.
\vspace{2ex}

\noindent {\bf Example 1:} \ We have mentioned earlier that the radial part of
a Brownian motion
$x_t$ in $X=\bR^n$ ($n\geq 2$), under the action of $K=SO(n)$, is a Bessel process
in a fixed ray $Y$ from the origin. We may take $Y$ to be the positive half of
$x_1$-axis, which is transversal to $K$ with boundary containing only the origin.
Then $M=\diag\{1,\,SO(n-1)\}$.
\vspace{2ex}

\noindent {\bf Example 2:} \ Let $X$ be the space of $n\times n$ real symmetric matrices
($n\geq 2$) with $K=SO(n)$ acting on $X$ by conjugation. The set $Y$
of all $n\times n$ diagonal matrices with
non-ascending diagonal elements is a submanifold of $X$ transversal to $SO(n)$,
its boundary consists of diagonal matrices with two identical diagonal elements,
and $M$ is the finite subgroup of $SO(n)$ consisting of diagonal matrices
with $\pm 1$ along diagonal.
The map $J$: $X\to Y$ maps a symmetric matrix to the diagonal matrix of its eigenvalues
in non-ascending order. Note that $X=GL(n,\bR)/SO(n)$, where $GL(n,\bR)$ is
the group of $n\times n$ real invertible matrices.
\vspace{2ex}

\noindent {\bf Example 3:} \ Let $Y$ be a manifold and $K$ be a Lie group with
a compact subgroup $M$, and let $X=Y\times(K/M)$ as a product manifold.
Then $K$ acts on $X$ as its natural action on $K/M$, $Y$ is transversal to $K$ and
$M$ is the isotropy subgroup of $K$ at all $y\in Y$. For example, $X=\bR^{n+m}=
Y\times K$ with $Y=\bR^n$, $K=\bR^m$ (additive group) and $M=\{0\}$.
\vspace{2ex}

\noindent {\bf Example 4:} \ Let $X=S^n$ be the $n$-dimensional sphere, regarded as
the unit sphere in $\bR^{n+1}$, under the natural action of $K=\diag\{1,SO(n)\}$.
The half circle $Y$ connecting two poles $(\pm 1,0,\ldots,0)$,
given by $(\cos t,\sin t,0,\ldots,0)$ for $0\leq t\leq\pi$, is transversal to $K$,
and $M=\diag\{1,1,\,SO(n-2)\}$.
\vspace{2ex}

\noindent {\bf Example 5:} \ Let $X=G/K$ be a symmetric space of noncompact type,
where $G$ is a semisimple Lie group of noncompact type and of a finite center,
and $K$ is a maximal compact subgroup. Using the standard notation and results
in \cite{helgason1}, let $\myg$ and $\myk$ be respectively the Lie algebras
of $G$ and $K$, let $\myp$ be an $\Ad(K)$-invariant subspace of $\myg$ complementary
to $\myk$, let $\mya$ be a maximal abelian subspace
of $\myp$, and let $\mya_+$ be a fixed (open) Weyl chamber ($\subset\mya$).
Then $Y=\exp(\overline{\mya_+})$ (overline denotes the closure) is a submanifold of $X$
transversal to the action of $K$ on $G/K$, its boundary is $\exp(\partial\mya_+)$,
where $\partial\mya_+$ is the boundary of $\mya_+$, and the isotropy subgroup $M$ of $K$ at
any $y\in Y^\circ=\exp(\mya_+)$ is the centralizer $M$ of $A$ in $K$.

\section{Nonhomogeneous \levy processes} \label{nonh}

Let $G$ be a Lie group. The convolution of two finite measures $\mu$ and $\nu$
on $G$ is the measure $\mu*\nu$ defined by $\mu*\nu(f)=\int f(xy)
\mu(dx)\nu(dy)$ for $f\in C_b(G)$. A family of probability measures
$\mu_t$, $t\in\bR_+$, on $G$ is called a convolution semigroup if $\mu_{t+s}=
\mu_t*\mu_s$ and $\mu_0=\delta_e$, the unit point mass at the identity element $e$
of $G$. It is called continuous if $\mu_t\to\mu_0$ weakly as $t\to 0$.
If $x_t$ is a \levy process in $G$ with transition semigroup $P_t$,
then $\mu_t=P_t(e,\cdot)$ is a continuous convolution semigroup. Conversely,
any continuous convolution semigroup is associated to
a \levy process in $G$ in this way (see \cite[chapter 1]{liao}).

Let $H$ be a compact subgroup of $G$. The convolution product may be extended to the
homogeneous space $G/H$ as follows. Let $\pi$: $G\to G/H$ be the natural projection.
A measurable map $S$: $G/H\to G$ is called a section map if $\pi\circ S=\id_{G/H}$,
the identity map on $G/H$. The convolution of two $H$-invariant finite measures
$\mu$ and $\nu$ on $G/H$ is the measure $\mu*\nu$ defined by $\mu*\nu(f)=\int
f(S(x)y)\mu(dx)\nu(dy)$, which is independent of the choice
of $S$ and is $H$-invariant. Moreover, the convolution product is associative and
hence the $n$-fold product $\mu_1*\mu_2*\cdots*\mu_n$ is well defined.
The convolution semigroup of $H$-invariant probability
measures on $G/H$ is defined as on $G$ but replacing $e$ by $o=eH$ (the origin of $G/H$).
The continuous convolution semigroups on $G/H$ are associated to
\levy processes as on $G$.

A process $x_t$ in a Lie group $G$ with rcll paths and an infinite life time
is called a nonhomogeneous \levy process
if for $s<t$, its increment $x_s^{-1}x_t$ is independent of process up to time $s$.
The distributions $\mu_{s,t}$ of the increments $x_s^{-1}x_t$, $s\leq t$, form
a two-parameter convolution semigroup in the sense that for $s<t<u$,
$\mu_{s,t}*\mu_{t,u}=\mu_{s,u}$ and $\mu_{t,t}=\delta_e$, which is continuous in
the sense that $\mu_{s,t}\to\mu_{s,s}$ weakly as $t\downarrow s$. In fact,
a nonhomogeneous \levy process in $G$ may be defined as a rcll process $x_t$ such that
for $0<t_1<t_2<\cdots<t_n$ and $f\in C_b(G^{n+1})$,
\begin{eqnarray}
&& E[f(x_{t_0},x_{t_1},x_{t_2}\ldots,x_{t_n})] = \int
f(x_0,x_0x_1,x_0x_1x_2,\ldots,x_0x_1\cdots x_n) \nonumber \\
&& \hspace{2in} \mu_0(dx_0)\mu_{0,t_1}
(dx_1)\mu_{t_1,t_2}(dx_2)\cdots\mu_{t_{n-1},t_n}(dx_n) \label{xtdist}
\end{eqnarray}
for a probability measure $\mu_0$ (distribution of $x_0$) and a continuous
two-parameter convolution semigroup $\mu_{s,t}$ on $G$.

Feinsilver \cite{feinsilver} obtained a martingale representation of such processes.
Let $\myg$ be the Lie algebra of $G$, whose elements are identified with left invariant
vector fields on $G$ as usual, and let $\xi_1,\ldots,\xi_n$ be a basis of $\myg$. Choose
local coordinates $\phi_1,\ldots,\phi_n\in C_c^\infty(G)$ to satisfy
$x=\exp(\sum_i\phi_i(x)\xi_i)$ for $x$ near $e$.
A covariance function $A$ is a continuous $n\times n$ symmetric matrix valued function such
that $A(0)=0$ and for $s<t$, $A(t)-A(s)$ is nonnegative definite. A \levy measure function
$\Pi(t,dx)$ is a measure valued function on $G$ such that $\Pi(0,\cdot)=0$,
$\Pi(t,\{e\})=0$ and for $f\in C_b^\infty(G)$ with $f(e)=\xi_if(e)=0$,
$\Pi(t,f)$ is finite and continuous
in $t$. Let $x_t$ be a nonhomogeneous \levy process in $G$ with $x_0=e$. Assume $x_t$
is stochastic continuous (that is, $x_t=x_{t-}$ almost surely for each fixed $t$).
Then by \cite{feinsilver}, there are unique $G$-valued
(non-random) continuous function $b_t$ with $b_0=e$, covariance function $A$ and
\levy measure function $\Pi$, such that $x_t=z_tb_t$ and for $f\in C_c^\infty(G)$,
\begin{eqnarray}
f(z_t) &-& \int_0^t\int_G\{f(z_sb_s\tau b_s^{-1})-f(z_s)-\sum_i\phi_i(\tau)
[\Ad(b_s)\xi_i]f(z_s)\}\Pi(ds,d\tau) \nonumber \\
&& \hspace{-0.3in}
- \int_0^t\frac{1}{2}\sum_{i,j}[\Ad(b_s)\xi_i][\Ad(b_s)\xi_j]f(z_s)\,dA_{ij}(s)
\label{mart}
\end{eqnarray}
is a martingale. Moreover, given $(b,A,\Pi)$ as above, there is a rcll process $x_t=z_tb_t$
in $G$ with $x_0=e$ such that (\ref{mart}) is a martingale for $f\in C_c^\infty(G)$.
Furthermore, such a process $x_t$ is unique in distribution and is a stochastic continuous
nonhomogeneous \levy process in $G$.

Note that because the exponential coordinates $\phi_i$ are used,
$\rho^{ij}_k= \xi_i\xi_j\phi_k(e)$ satisfies $\rho^{ij}_k=-\rho^{ji}_k$, so they
will not appear in (\ref{mart}) as in \cite{feinsilver}. Note also that
the integrand of $\Pi$-integral in (\ref{mart}) is the remainder of
a first order Taylor expansion and hence is integrable.

Now consider the homogeneous space $G/H$. A point $b\in G/H$ or a subset $B$
of $G/H$ is called $H$-invariant if $hb=b$ or $hB=B$ for all $h\in H$.
For $x\in G/H$,
$H$-invariant $b$ and $B$, $xb=S(x)b\in G/H$ and $xB=S(x)B\subset G/H$ are
well defined because they are independent of choice for the section map $S$.
Note that $g\in G$ with $go$ $H$-invariant is characterized by $g^{-1}Hg\subset H$, and
hence by $g^{-1}Hg=H$. Therefore, the set of $H$-invariant points in $G/H$ is
the natural projection of a closed subgroup of $G$ containing $H$ as a normal subgroup,
and hence has a natural group structure with product $b_1b_2=S(b_1)b_2$ and
inverse $b^{-1}=S(b)^{-1}o$ (independent of $S$).
In general, the product $x_1x_2\cdots x_{n-1}x_n=S(x_1)S(x_2)\cdots S(x_{n-1})x_n$
is not well defined because it depends on the choice of $S$.
However, if $x_2,\ldots,x_n$ are independent random variables
with $H$-invariant distributions, then the distribution of the product $x_1x_2\cdots
x_n$ and also that of the sequence $y_i=x_1x_2\cdots x_i$ for $i=1,2,\ldots$
are independent of $S$, and hence such a product or sequence
is meaningful in the sense of distribution. 

Note that for $H$-invariant finite measures $\mu$ and $\nu$ on $G/H$, an integral like
\[\int f(xy,xyz)\mu(dy)\nu(dz) = \int
f(S(x)y,\,S(x)S(y)z)\mu(dy)\nu(dz)\]
is well defined (independent of choice of section map $S$). So is
$\int f(xbyb^{-1})\mu(dy)$ if $b$ is an $H$-invariant point.

A process $x_t$ in $G/H$ with rcll paths and an infinite life time
will be called a nonhomogeneous \levy process
if there is a continuous two-parameter convolution semigroups $\mu_{s,t}$
of $H$-invariant probability measures on $G/H$ such that (\ref{xtdist}) holds.
Then for $s<t$, $x_s^{-1}x_t=S(x_s)^{-1}x_t$
has distribution $\mu_{s,t}$ (independent of choice for section map $S$)
and is independent of the process up to time $s$.

Because $H$ is compact, there is a subspace $\myp$ of $\myg$
that is complementary to the Lie algebra $\myh$ of $H$ and is $\Ad(H)$-invariant
in the sense that $\Ad(h)\myp=\myp$ for $h\in H$. Choose a basis $\xi_1,
\ldots,\xi_m$ of $\myp$ and local coordinates $\phi_1,\ldots,\phi_m\in C_c^\infty(G/H)$
around $o$ on $G/H$ such that $x=\exp(\sum_{i=1}^m \phi_i(x)\xi_i)o$ for $x$ near $o$.
Then
\begin{equation}
\forall h\in H,\ \ \ \ \sum_{i=1}^m \phi_i\Ad(h)\xi_i = \sum_{i=1}^m(\phi_i\circ h)\xi_i
\label{xiAdh}
\end{equation}
near $o$. The functions $\phi_i$ may be suitably extended so that (\ref{xiAdh}) holds
globally on $G/H$.

Any $\xi\in\myg$ is a left invariant vector field on $G$.
If $\xi$ is $\Ad(H)$-invariant, it may also be regarded as a vector field on $G/H$
given by $\xi f(x)=\frac{d}{dt}f(xe^{t\xi}o)\mid_{t=0}$ for $f\in C^\infty(G/H)$
and $x\in G/H$ (note that $e^{t\xi}o$ is $H$-invariant), which is
$G$-invariant in the sense that $\xi(f\circ g)=(\xi f)\circ g$ for $g\in G$.
In fact, any $G$-invariant vector field on $G/H$ is given by an $\Ad(H)$-invariant
$\xi\in\myg$. Note that if $\xi\in\myg$ is $\Ad(H)$-invariant and $b\in G/H$
is $H$-invariant, then $\Ad(b)\xi=\Ad(S(b))\xi$ is $\Ad(H)$-invariant and is independent
of section map $S$. By (\ref{xiAdh}), for any $H$-invariant measure $\mu$ on $G/H$,
$\int\mu(dx)\sum_i\phi_i(x)\xi_i$ is $\Ad(H)$-invariant, and so is
$\int\mu(dx)\sum_i\phi_i(x)\Ad(b)\xi_i$.

Let $\xi,\eta\in\myg$. With a choice of section map $S$, $\xi\eta$ may be regarded
as a second order differential operator
on $G/H$ by setting $\xi\eta f(x)=\frac{\partial^2}{\partial t\,\partial s}
f(S(x)e^{t\xi}e^{s\eta}o)\mid_{t=s=0}$.
As in \cite{feinsilver}, it can be shown that
$\xi_i\xi_j f(x) = \frac{\partial^2}{\partial t\,\partial s}
f(S(x)e^{t\xi_i+s\xi_j}o)\mid_{t=s=0}+\sum_{k=1}^m\rho_k^{ij}\xi_kf(x)$
with $\rho_k^{ij}=-\rho_k^{ji}$. Thus, if $a_{ij}$ is a symmetric matrix, then
\[\sum_{i,j=1}^ma_{ij}\xi_i\xi_jf(x) =
\sum_{i,j=1}^ma_{ij}\frac{\partial^2}{\partial t_i\,\partial t_j}
f(S(x)e^{\sum_{p=1}^m t_p\xi_p}o)\mid_{t_1=\cdots=t_m=0}.\]
The matrix $a_{ij}$ is called $\Ad(H)$-invariant
if $a_{ij}=\sum_{p,q}a_{pq}[\Ad(h)]_{pi}[\Ad(h)]_{qj}$ for $h\in H$, where
$[\Ad(h)]_{ij}$ is the matrix representing $\Ad(h)$, that is,
$\Ad(h)\xi_j=\sum_i[\Ad(h)]_{ij}\xi_j$. Then the operator $\sum_{i,j}a_{ij}
\xi_i\xi_j$ is independent of section map $S$ and is $G$-invariant.
In fact, any second order $G$-invariant differential operator
on $G/H$ is such an operator plus a $G$-invariant vector field. Note that if $b\in G/H$
is $H$-invariant, then $\sum_{i,j}a_{ij}[\Ad(b)\xi_i][\Ad(b)\xi_j]=\sum_{i,j}a_{ij}
[\Ad(S(b))\xi_i][\Ad(S(b))\xi_j]$ is a $G$-invariant operator on $G/H$
(independent of $S$).

A covariance function $A$ and a \levy measure function $\Pi$ on $G/H$ are defined as on $G$
with the additional requirements that $A(t)$ is $\Ad(H)$-invariant and $\Pi(t,\cdot)$
is $H$-invariant. By the preceding discussion,
the expression in (\ref{mart}) is meaningful on $G/H$ and is independent of the choice
of section map $S$ in $\Ad(b_s)=\Ad(S(b_s))$.
The following result is an extension of Feinsilver's martingale
representation to nonhomogeneous \levy processes in $G/H$.

\begin{theo} \label{th2}
Let $x_t$ be a stochastic continuous nonhomogeneous \levy process in $G/H$ with $x_0=o$.
Then there is a unique triple $(b,A,\Pi)$ of a continuous $H$-invariant function
$b_t$ with $b_0=o$, a covariance function $A$ and a \levy measure function $\Pi$
on $G/H$ such that $x_t=z_tb_t$ and (\ref{mart}) is a martingale for $f\in C_c^\infty(G/H)$.
Moreover, given $(b,A,\Pi)$ as above, there is a rcll process $x_t=z_tb_t$
in $G/H$ with $x_0=o$ and represented by $(b,A,\Pi)$. Furthermore, such $x_t$
is unique in distribution and is a stochastic continuous nonhomogeneous process in $G/H$.
\end{theo}

\noindent {\bf Proof} \ Fix $T>0$. Let $0=t_0<t_1<\cdots<t_n\leq T$ be a partition
of $[0,\,T]$ with $t_{i+1}-t_i=1/n$ and let $\mu_{ni}=\mu_{t_{i-1},t_i}$
for $1\leq i\leq n$. By the stochastic continuity of $x_t$, $\int \phi_id\mu_{ni}$ is
uniformly small in $i$ as $n\to\infty$, and $b_{ni}\in G/H$ given by
$x_i(b_{ni})=\int \phi_id\mu_{ni}$ is well defined and is $H$-invariant because
the $H$-invariance of $\mu_{ni}$.

Define an $H$-invariant function $b_n(t)$ in $G/H$
by $b_n(t)=b_{n1}b_{n2}\cdots b_{n[nt]}$ for $0<t\leq T$ and $b_n(0)=o$,
where $[nt]$ is the integer part of $nt$,
a measure function $\Pi_n$ by $\Pi_n(t,\cdot)=\sum_{i=1}^{[nt]}\mu_{ni}$ for $0<t\leq T$
and $\Pi(0,\cdot)=0$, and a matrix valued function $A_n(t,U)$
of time $t$ and measurable $U\subset G/H$
by $A_n(0,U)=0$ and for $0<t\leq T$,
\[[A_n(t,U)]_{ij} = \sum_{k=1}^{[nt]}\int_U[\phi_i(x)-\phi_i(b_{nk})]
[\phi_j(x)-\phi_j(b_{nk})]\mu_{nk}(dx).\]

Let $x_{ni}$, $1\leq i\leq n$, be independent random variables in $G/H$
with distributions $\mu_{ni}$. Define $x_n(t)=x_{n1}x_{n2}\cdots x_{n[nt]}$
and $z_n(t)=x_n(t)b_n(t)^{-1}$. Then as $n\to\infty$, the process $x_n(t)$ converges
in distribution to $x_t$, and it can be shown that
$b_n(t)$ converges uniformly to a continuous $H$-invariant function
$b_t$, $\Pi_n$ converges to a \levy measure function $\Pi$ in the sense that
for any $f\in C_b(G/H)$ vanishing near $o$, $\Pi_n(t,f)\to\Pi(t,f)$ uniformly
for $0\leq t\leq T$, for any $H$-invariant neighborhood $U$ of $o$ which is
a continuity set of $\Pi(T,\cdot)$, $[A_n(t,U)]_{ij}\to A_{ij}(t)+\int_U
\phi_i(x)\phi_j(x)\Pi(t,dx)$ uniformly for $0\leq t\leq T$, for some covariance function
$A_{ij}(t)$, and $z_n(t)$ converges in distribution to a stochastic continuous
process $z_t$ for which (\ref{mart}) is a martingale. Thus, $(b,A,\Pi)$ are
the parameters in the representation of $x_t$.
These statements and the rest of theorem may be proved by essentially
repeating the proof in \cite{feinsilver} with suitable changes, such as
properly interpreting the product on $G/H$ as discussed earlier and
using $H$-invariant sets on $G/H$ for various neighborhoods used
in \cite{feinsilver}. \ $\Box$
\vspace{2ex}

Let $x_t$ be a stochastic continuous nonhomogeneous \levy process in $G/H$.
Fix $T>0$. We now show that for $f\in C_b(G/H)$ with support supp$(f)$ not containing $o$,
\begin{equation}
\Pi(T,f) = E[\sum_{0<t\leq T}f(x_{t-}^{-1}x_t)], \label{PiTf}
\end{equation}
which is independent of choice for section map $S$ to represent
$x_{t-}^{-1}x_t=S(x_{t-})^{-1}x_t$, where the summation $\sum_{0<t\leq T}
f(x_{t-}^{-1}x_t)$ has only finitely many nonzero terms almost surely
because $o\not\in$supp$(f)$.
Take a partition $0=t_0<t_1<\cdots<t_n\leq T$ of $[0,\,T]$ as in the proof
of Theorem~\ref{th2}. Let $\psi_j$, $j=1,2,\ldots$, be a partition of unity on $G/H$,
that is, $0\leq \psi_j\in C_c(G/H)$ and $\sum_j\psi_j=1$ is a locally finite sum (that is,
any point of $G/H$ has a neighborhood on which only finitely many $\psi_j$ are nonzero).
We may assume that for each $j$, there is a section map $S_j$ continuous
on supp$(\psi_j)$. Then $E[\sum_i\psi_j(x_{t_{i-1}})
f(S_j(x_{t_{i-1}})^{-1}x_{t_i})]\to E[\sum_{0<t\leq T}\psi_j(x_{t-})
f(S_j(x_{t-})^{-1}x_t)]$ as $n\to\infty$. This holds for any section
map $S_j$ continuous on supp$(\psi_j)$, and by (\ref{xtdist}) and the $H$-invariance
of $\mu_{s,t}$, these expressions are independent of choice of such an $S_j$.
It then follows that $S_j(x)$ in these expressions may be replaced by $S(x)$
for any section map $S$.
By the proof of Theorem~\ref{th2} and the local finiteness of $\sum_j\psi_j=1$,
\begin{eqnarray*}
&& \Pi(T,f) \leftarrow \sum_{i=1}^{[nt]}\mu_{ni}(f) = \sum_j\sum_i
E[\psi_j(x_{t_{i-1}})f(S(x_{t_{i-1}})^{-1}x_{t_i})] \\
&=& \sum_j\sum_i
E[\psi_j(x_{t_{i-1}})f(S_j(x_{t_{i-1}})^{-1}x_{t_i})]
\to \sum_j E[\sum_{0<t\leq T}\psi_j(x_{t-})f(S_j(x_{t-})^{-1}x_t)] \\
&=& \sum_j E[\sum_{0<t\leq T}\psi_j(x_{t-})f(S(x_{t-})^{-1}x_t)] =
E[\sum_{0<t\leq T}f(S(x_{t-})^{-1}x_t)]. \ \ \ \ \Box
\end{eqnarray*}

By (\ref{PiTf}), $\Pi(T,B)$ is the expected numbers of jumps
$x_{t-}^{-1}x_t=S(x_{t-})^{-1}x_t$ in $B\subset G/H$ for $0<t\leq T$ (independent of $S$).
Thus, $x_t$ is continuous if and only if $\Pi=0$.
It can be shown that the pairs $(t,\,x_{t-}^{-1}x_t)$
with $x_{t-}^{-1}x_t\neq o$, $0<t\leq T$, form a Poisson random measure $N$
on $[0,\,T]\times(G/H)$ with intensity measure $\Pi(dt,dx)$.

Note that if $\Pi(t,\vert f\vert)<\infty$ for any $f\in C_b(G/H)$ with $f(o)=0$
and $t>0$, in particular, if $\Pi(t,\cdot)$ is a finite measure, then
by suitably changing the drift $b_t$, the term $\sum_ix_i[\Ad(b_s)\xi_i]f(z_s)$
in the integrand of $\Pi$-integral in (\ref{mart}) may be dropped. To prove this,
apply It\^{o}'s formula to $f(z_tu_t)$, where $u_t$ is the $H$-invariant
function in $G/H$ determined by the ordinary differential equation $du_t=\int_0^t\int_{G/H}
\Pi(ds,d\tau)\sum_i\phi_i(\tau)[\Ad(b_s)\xi_i](u_s)$ and $u_0=o$.

The \levy measure function $\Pi$ is clearly independent of the choice for
basis $\xi_1,\ldots, \xi_m$ of $\myp$. Although the covariance function $A$ depend
on the basis, its uniqueness under a given basis implies that the differential
operator $V(t)=(1/2)\sum_{i,j=1}^mA_{ij}(t)\xi_i\xi_j$ on $G/H$ is independent of
the basis. The operator $V(t)$ will be called the covariance operator of $x_t$,
which together with $b_t$ and $\Pi(t,\cdot)$ determines the distribution of
the process $x_t$ completely.

In general, a nonhomogeneous \levy process $x_t$ in $G/H$ may have a fixed jump,
that is, $P(x_{t-}\neq x_t)>0$ for some $t>0$. Suppose $x_t$ have only finitely
many fixed jumps at times $t_1<t_2<\cdots<t_k$. As jumps $x_{t-}^{-1}x_t$
for $t=t_j$ are independent of the rest of process $x_t$, they may be easily removed to
obtain a stochastic continuous nonhomogeneous \levy process $x_t'$. The triple
$(b,A,\Pi)$ (or $(b,V,\Pi)$ with $V$ in place of $A$) of $x_t'$ will also be called
the drift, covariance function (or covariance operator) and \levy measure function
of $x_t$, which together with the distributions of fixed jumps determine
the distribution of process $x_t$ completely.

A consequence of Theorem~\ref{th2} is that a nonhomogeneous
\levy process $x_t$ in $G/H$ is the projection of such a process $g_t$ in $G$
which is $H$-conjugate invariant in the sense that for any $h\in H$,
the process $c_h(g_t)$ has the same distribution as $g_t$, where $c_h$:
$G\ni g\mapsto hgh^{-1}\in G$ is the conjugation map. However, such $g_t$ is
not unique in distribution. A similar result for (homogeneous) \levy processes
are obtained in \cite[chapter 2]{liao}.

\begin{coro} \label{co1}
Let $x_t$ be a stochastic continuous nonhomogeneous \levy process in $G/H$ with
$x_0=o$. Then there is a stochastic continuous $H$-conjugate invariant \levy process
$g_t$ in $G$ with $g_0=e$ such that the two processes $x_t$ and $g_to$ are
identical in distribution.
\end{coro}

\noindent {\bf Proof} \ Let $(b,A,\Pi)$ be the representation of $x_t$ in Theorem~\ref{th2}.
There is a continuous $G$-valued function $b_t'$ in $G$ with $b_0'=e$ and
$b_t=b_t'o$. Let $\hat{b}_t=\int_H hb_t'h^{-1}dh$, where $dh$ is the normalized Haar
measure on $H$. Then $\hat{b}_t$ is a continuous $H$-conjugate invariant function
in $G$ with $\hat{b}_0=e$ and $b_t=\hat{b}_to$. The basis $\xi_1,\ldots,\xi_m$ of $\myp$
may be extended to be a basis $\xi_1,\ldots,\xi_n$ of $\myg$ such that $\xi_{m+1},
\ldots,\xi_n$ form a basis of $\myh$.
The covariance function $A_{ij}(t)$ on $G/H$ may be regarded as covariance function on $G$
with $A_{ij}(t)=0$ if either $i>m$ or $j>m$.
Let $S$ be a section map satisfying $S(x)=\exp[\sum_{i=1}^m
\phi_i(x)\xi_i]$ for $x$ near $o$, and let $\hat{\Pi}(t,\cdot)$ be defined by
$\hat{\Pi}(t,f)=\int_H f(hS(x)h^{-1})dh\Pi(t,dx)$ for $f\in C_b(G)$. Then $\hat{\Pi}(t,\cdot)$
is a $H$-conjugate invariant \levy measure function on $G$ (that is, $c_h
\hat{\Pi}(t,\cdot)=\hat{\Pi}(t,\cdot)$ for $h\in H$). The nonhomogeneous
\levy process $g_t$ in $G$ with $g_0=e$ and representation $(\hat{b},A,\hat{\Pi})$ satisfies
$x_t=g_to$ in distribution because $(\hat{b},A,\hat{\Pi})$ on $G$
project to $(b,A,\Pi)$ on $G/H$, and is $H$-conjugate invariant because both $\hat{b}_t$
and $\hat{\Pi}(t,\cdot)$ are $H$-conjugate invariant, and $A(t)$ is $\Ad(H)$-invariant. \ $\Box$


\section{Angular part} \label{angular}

Recall that $x_t$ is a Markov process in a manifold $X$ invariant under the action of
a Lie group $K$, $Y$ is a submanifold transversal to $K$ with interior $Y^\circ$, $M$ is
the compact isotropy subgroup at every point of $Y^\circ$, and $X^\circ$ is
the union of $K$-orbits through $Y^\circ$ with $X^\circ=Y^\circ\times(K/M)$.

The exit time $\zeta$ of process $x_t$ from $X^\circ$ is the stopping time when $x_t$
together with its left limit first leaves $X^\circ$ or reaches its life time $\xi$.
More precisely, it is defined by
\begin{equation}
\zeta = \inf\{t>0;\ \ \ \ x_t\not\in X^\circ,\ x_{t-}\not\in X^\circ\ {\rm or}
\ t\geq\xi\}, \label{zeta}
\end{equation}
with $\inf$ of an empty set defined to be $\infty$. The exit time of $y_t$ from $Y^\circ$
is also denoted by $\zeta$. Fix $T>0$. Because the process $x_t$ has rcll paths,
it may be regarded as a random variable in the space $D_T(X)$
of rcll maps: $[0,\,T]\to X^\circ$, equipped with Skorohod topology.
Let $P_x$ be the distribution on $D_T(X)$ associated to the process $x_t$ starting
at $x\in X$. Its total mass may be less than $1$ because $x_t$ may have a finite life time.

For $y\in Y^\circ$ and $z\in K/M$, $zy=S(z)y\in X^\circ$ is well defined and is independent
of choice of section map $S$: $K/M\to K$. Let $x_t=z_ty_t$ be the decomposition of
the process $x_t$ with $x_0\in X^\circ$ and $t<\zeta$. Then $y_t$ is the radial part
as defined before, and $z_t$ is a process in $K/M$ with rcll paths and will be called
the angular part of $x_t$.

Recall $J$ is the projection map $X\ni x\mapsto y\in Y$.
Let $J_2$ be the projection map $X^\circ\ni x\mapsto z\in K/M$
associated to the decomposition $x=zy$. We will also use $J$ and $J_2$
to denote the maps $J$: $D_T(X)\ni x(\cdot)\mapsto y(\cdot)\in D_T(Y)$
and $J_2$: $D_T(X^\circ)\ni x(\cdot)\mapsto z(\cdot)\in D_T(K/M)$ respectively
given by the decomposition $x(\cdot)=z(\cdot)y(\cdot)$.

Let $\myF_{0,T}^Y=\sigma\{y_t$; $0\leq t\leq T\}$ be the $\sigma$-algebra generated by
the radial process $y_t$ for $0\leq t\leq T$, which may be regarded as
a $\sigma$-algebra on $D_T(Y)$ and induces the $\sigma$-algebra $J^{-1}(\myF_{0,T}^Y)$
on $D_T(X)$. By the existence of regular conditional distributions (see
for example \cite[chapter 5]{kallenberg}), there is a probability kernel $R_z^{y(\cdot)}$
from $D_T(Y^\circ)\times(K/M)$ to $D_T(K/M)$ such that
for any $x\in X^\circ$ and measurable $F\subset D_T(K/M)$,
\begin{equation}
R_{J_2(x)}^{J[x(\cdot)]}(F) = P_x[J_2^{-1}(F)\mid J^{-1}(\myF_{0,T}^Y)]\ \ \ \ \mbox{
for $P_x$-almost all $x(\cdot)$ in $[\zeta>T]\subset D_T(X^\circ)$}. \label{Qy}
\end{equation}
The probability measure $R_z^{y(\cdot)}$ is the conditional distribution of
the angular process $z_t$ given a radial path $y(\cdot)$ in $D_T(Y^\circ)$ and $z_0=z$.

\begin{theo} \label{th3}
Fix $T>0$. Almost surely on $[\zeta>T]$, given a radial path $y_t$ for $0\leq t\leq T$,
the conditioned angular process $z_t$ is a nonhomogeneous \levy process in $K/M$.
More precisely, this means that for $y\in Y^\circ$, $z\in K/M$, and $JP_y$-almost
all $y(\cdot)$ in $[\zeta>T]\subset D_T(Y^\circ)$, the angular process $z_t$ is
a nonhomogeneous \levy process under $R_z^{y(\cdot)}$.
\end{theo}

\noindent {\bf Proof} \ For $x\in X^\circ$, let $\tilde{P}_t(x,B)=
P_x\{[x_t\in B]\cap[\zeta>t]\}$ for measurable $B\subset X^\circ$. By the simple
Markov property of $x_t$, it is easy to show that $\tilde{P}_t$ is
the transition semigroup of the Markov process $x_t$ for $t<\zeta$ and
it is $K$-invariant. Similarly,
let $\tilde{Q}_t$ be the transition semigroup of $y_t$ for $t<\zeta$.
Then $\tilde{P}_t(x,\cdot)$ and $\tilde{Q}_t(y,\cdot)$ are
respectively sub-probability kernels from $X^\circ$ to $X^\circ=Y^\circ
\times(K/M)$ and from $Y^\circ$ to $Y^\circ$. By the existence of
a regular conditional distribution, there is a probability kernel $R_t(y,y_1,\cdot)$
from $(Y^\circ)^2$ to $K/M$ such that for $y\in Y^\circ$,
$\tilde{P}_t(y,\,dy_1\times dz_1)=\tilde{Q}_t(y,dy_1)R_t(y,y_1,dz_1)$.
The $K$-invariance of $\tilde{P}_t$ implies that the measure $R_t(y,y_1,\cdot)$
is $M$-invariant for $Q_t(y,\cdot)$-almost all $y_1$. Modifying $R_t$ on
an exceptional set
of zero $Q_t(y,\cdot)$-measure, we may assume $R_t(y,y_1,\cdot)$ is $M$-invariant
for all $y,y_1\in Y^\circ$. Therefore,
for $z\in K/M$, it is meaningful to write $R_t(y,y_1,z^{-1}dz_1)=R_t(y,y_1,S(z)^{-1}dz_1)$
because it is independent of choice of section map $S$. We then have
\begin{equation}
\forall y\in Y^\circ\ {\rm and}\ z\in K/M,
\ \ \ \ \tilde{P}_t(zy,\,dy_1\times dz_1) = \tilde{Q}_t(y,dy_1)R_t(y,y_1,z^{-1}dz_1).
\label{PtQt}
\end{equation}

For $0<s_1<s_2<\cdots<s_k<\infty$, $y\in Y^\circ$, $z\in K/M$, $h\in C_b(Y^k)$
and $f\in C_b((K/M)^k)$,
\begin{eqnarray*}
&& E_{zy}[h(y_{s_1},\ldots,y_{s_k})f(z_{s_1},\ldots,z_{s_k});\,\zeta>s_k] \nonumber \\
&=& \int\int \tilde{P}_{s_1}(zy,\,dy_1\times dz_1)
\tilde{P}_{s_2-s_1}(z_1y_1,\,dy_2\times dz_2) \cdots
\tilde{P}_{s_k-s_{k-1}}(z_{k-1}y_{k-1},\,dy_k\times dz_k) \nonumber \\
&& \ \ h(y_1,y_2,\ldots,y_k)f(z_1,z_2,\ldots,z_k) \nonumber \\
&=& \int \tilde{Q}_{s_1}(y,dy_1)\tilde{Q}_{s_2-s_1}(y_1,dy_2)\cdots
\tilde{Q}_{s_k-s_{k-1}}(y_{k-1},dy_k)h(y_1,y_2,\ldots,y_k) \nonumber \\
&& \hspace{-0.2in} \int R_{s_1}(y,y_1,dz_1)R_{s_2-s_1}(y_1,y_2,dz_2)\cdots R_{s_k-s_{k-1}}
(y_{k-1},y_k,dz_k)f(zz_1,zz_1z_2,\ldots,zz_1\cdots z_k) \nonumber \\
&=& E_y[h(y_{s_1},y_{s_2},\ldots,y_{s_k})\int R_{s_1}(y,y_{s_1},dz_1)
R_{s_2-s_1}(y_{s_1},y_{s_2},dz_2)\cdots \nonumber \\
&& \hspace{1.7in} R_{s_k-s_{k-1}}(y_{s_{k-1}},y_{s_k},dz_k)
f(zz_1,zz_1z_2,\ldots,zz_1\cdots z_k);\,\zeta>s_k].
\end{eqnarray*}
This implies that on $[\zeta>s_k]$,
\begin{eqnarray}
&& E_{zy}[f(z_{s_1},\ldots,z_{s_k})\mid y_{s_1},\ldots,y_{s_k}] =
\int R_{s_1}(y,y_{s_1},dz_1)R_{s_2-s_1}(y_{s_1},y_{s_2},dz_2)\cdots \nonumber \\
&& \hspace{1.5in} R_{s_k-s_{k-1}}(y_{s_{k-1}},y_{s_k},dz_k)f(zz_1,zz_1z_2,
\ldots,zz_1\cdots z_k)]. \label{Ehf}
\end{eqnarray}

Let $\Gamma$ be the set of dyadic numbers $i/2^m$ for integers $i\geq 0$ and $m>0$.
For the moment, assume $T\in\Gamma$.
For $s,t\in\Gamma$ with $s<t\leq T$, let $s=s_1<s_2<\cdots<s_k$ be a partition
of $[0,\,T]$ spaced by $1/2^m$ with $s=s_i$ and $t=s_j$, and let
\[\mu_{s,t}^m = R_{s_{i+1}-s_i}(y_{s_i},y_{s_{i+1}},\cdot)* R_{s_{i+2}-s_{i+1}}
(y_{s_{i+1}},y_{s_{i+2}},\cdot)
*\cdots * R_{s_j-s_{j-1}}(y_{s_{j-1}},y_{s_j},\cdot).\]
By (\ref{Ehf}), $M$-invariance of $P_t(y,y_1,\cdot)$ and the measurability of $\mu_{s,t}^m$
in $y_{s_i},\ldots,y_{s_j}$,
\begin{equation}
\mu_{s,t}^m(f) = E_{zy}[f(z_s^{-1}z_t)\mid y_{s_1},\ldots,y_{s_k}] =
E_{zy}[f(z_s^{-1}z_t)\mid y_{s_i},\ldots,y_{s_j}]\ \ \mbox{on $[\zeta>T]$} \label{mustm}
\end{equation}
for $f\in C_b(K/M)$,
which is independent of the choice for section map $S$ to represent
$z_s^{-1}z_t=S(z_s)^{-1}z_t$. By the right continuity of $y_t$,
as $m\to\infty$, $\sigma(y_{s_1},\ldots,y_{s_k}\}\uparrow\myF_{0,T}^Y$ and
$\sigma\{y_{s_i},\ldots,y_{s_j}\}\uparrow\myF_{s,t}^Y$,
it follows that as $m\to\infty$, almost surely, $\mu_{s,t}^m\to\mu_{s,t}$ weakly
for some $M$-invariant probability measure $\mu_{s,t}$ on $K/M$ such that
\begin{equation}
\forall f\in C_b(K/M),
\ \ \ \ \mu_{s,t}(f) = E_{zy}[f(z_s^{-1}z_t)\mid\myF_{0,T}^Y] =
E_{zy}[f(z_s^{-1}z_t)\mid\myF_{s,t}^Y]\ \ \mbox{on $[\zeta>T]$}. \label{must}
\end{equation}
Note that $\mu_{s,t}$ is an $\myF_{s,t}^Y$-measurable random measure independent
of starting point $zy$. Because $\Gamma$ is countable, the exception set of probability zero
in the above almost sure convergence may be chosen simultaneously for all $s<t$ in
$\Gamma$. Moreover, for $t_1<t_2<\cdots<t_n$ of $[0,\,T]$ in $\Gamma$,
it can be shown from (\ref{Ehf})
and by choosing a partition $s_1<s_2<\cdots<s_k$ from $\Gamma$ containing all
$t_i$ that almost surely on $[\zeta>T]$, for $f\in C_b((K/M)^n)$,
\begin{eqnarray*}
&& E_{zy}[f(z_{t_1},\ldots,z_{t_n})\mid\myF_{0,T}^Y] = \lim_{m\to\infty}
E_{zy}[f(z_{t_1},\ldots,z_{t_n})\mid y_{s_1},\ldots,y_{s_k}] \\
&=& \lim_{m\to\infty}\int f(zz_1,zz_1z_2,\ldots,zz_1\cdots z_n)\mu_{0,t_1}^m(dz_1)
\mu_{t_1,t_2}^m(dz_2)\cdots\mu_{t_{n-1},t_n}^m(dz_n).
\end{eqnarray*}
This implies that almost surely on $[\zeta>T]$, for $0\leq t_1<\cdots<t_n\leq T$ in $\Gamma$,
\begin{eqnarray}
&& E_{zy}[f(z_{t_1},\ldots,z_{t_n})\mid\myF_{0,T}^Y] \nonumber \\
&=& \int f(zz_1,zz_1z_2,\ldots,zz_1\cdots z_n)\mu_{0,t_1}(dz_1)
\mu_{t_2,t_1}(dz_2)\cdots\mu_{t_{n-1},t_n}(dz_n). \label{Emu}
\end{eqnarray}
In particular, $\mu_{s,t}$ for $s<t$ in $\Gamma$ form a two-parameter convolution semigroup
on $K/M$. To prove that the conditioned process $z_t$ is a nonhomogeneous \levy processes
in $K/M$, it remains to extend $\mu_{s,t}$ in (\ref{must}) to all real $s<t\leq T$ and
prove (\ref{Emu}) for real $0\leq t_1<\cdots<t_n\leq T$.

Let $f\in C_b(K/M)$. By the right continuity of $z_t$ and (\ref{must}),
for $s\in\Gamma$ and real $t$
with $s<t<T$, $\mu_{s,t}$ may be defined as the weak limit
of $\mu_{s,t_n}=P[S(z_s)^{-1}z_{t_n}\in\cdot\mid\myF_{0,T}^Y]$ as $\Gamma\ni t_n\downarrow t$,
which is independent of section map $S$.
For a real $s$, choose $\Gamma\ni s_n\downarrow s$. We can show
$\mu_{s_n,t}(f) = E[f(S(z_{s_n})^{-1}z_t)\mid\myF_{0,T}^Y] \to
E[f(S(z_s)^{-1}z_t)\mid\myF_{0,T}^Y]$ by using a partition of unity
$\{\psi_j\}$ and section map $S_j$ continuous on supp$(\psi_j)$ as in the
proof of (\ref{PiTf}) in section~\ref{nonh}.
We can then define $\mu_{s,t}(f)=\lim_n\mu_{s_n,t}(f)$.
Note that no additional exceptional set is produced in taking these limits.
It is easy to see that $\mu_{s,t}$ for real $0<s<t$ form
a two-parameter convolution semigroup on $K/M$ for which (\ref{must})
and (\ref{Emu}) hold. To show the conclusions are valid for any real $T>0$, we may choose
$T_m\in\Gamma$ with $T_m\downarrow T$. \ $\Box$
\vspace{2ex}

Because the natural action of $M$ on $K/M$ fixes $o$, it induces an action on
the tangent space $T_o(K/M)$ at $o$.
The homogeneous space $K/M$ will be called irreducible if the action of $M$ on
$T_o(K/M)$ is irreducible (that is, it has no nontrivial invariant subspace).
Among the examples mentioned in section~\ref{radial}, $K/M$ is irreducible in
Examples 1 and 4, and in Example~3 if it is chosen to be so, and in Example~5
if the symmetric space $G/K$ is of rank $1$ (see \cite{helgason2}).

If $K/M$ is irreducible, then, up to a constant multiple, there is a unique $M$-invariant
inner product on $T_o(K/M)$ (see for example Appendix~5 in \cite{kn}), and hence, there is
a unique second order $K$-invariant differential operator on $K/M$.
By choosing a $K$-invariant Riemannian metric on $K/M$,
which is unique up to a constant factor, any such operator is a multiple of
the Laplace-Beltrami operator $\Delta_{K/M}$ on $K/M$. The following result
is an extension of Galmarino's result mentioned in section~\ref{intro}.

\begin{theo} \label{th4}
Assume $K/M$ is irreducible. If $x_t$ is a continuous $K$-invariant Markov process in $X$
with radial part $y_t$ in $Y$,
then there are a Brownian motion $B(t)$ in $K/M$ under a $K$-invariant
Riemannian metric, independent of process $x_t$, and a real continuous non-decreasing
process $a_t$, with $a_0=0$ and $a_t-a_s$ $\myF_{s,t}^Y$-measurable
for $s<t$, such that the two processes $x_t$ and $B(a_t)y_t$, $t<\zeta$,
are identical in distribution.
\end{theo}

\noindent {\bf Proof} \ By Theorem~\ref{th3}, given a radial path $y_t$ for
$0\leq t\leq T$ with $[\zeta>T]$, the conditioned angular process $z_t$ is
a continuous nonhomogeneous \levy process in $K/M$. Let $(b,A,\Pi)$ be its representation
in Theorem~\ref{th2}. By the irreducibility of $K/M$,
there is no nonzero $M$-invariant tangent vector at $o$,
and hence, there is no $M$-invariant point near $o$ except $o$.
This implies that the drift $b_t=o$. Because $z_t$ is continuous,
$\Pi=0$. Because $\sum_{i,j}A_{ij}(t)\xi_i\xi_j$ is a $K$-invariant second
order differential operator on $K/M$, it must be equal to $a_t\Delta_{K/M}$ for
a continuous non-decreasing function $a_t$ with $a_0=0$. By (\ref{must}), $\mu_{s,t}$
is $\myF_{s,t}^Y$-measurable. From the construction of $A(t)$ from $\mu_{s,t}$
in the proof of Theorem~\ref{th2},
it is seen that $a_t-a_s$ is $\myF_{s,t}^Y$-measurable.
Let $B(t)$ be a Brownina motion in $K/M$ independent of process $x_t$.
It is enough to show that the conditioned
process $z_t$ is equal to the time-changed Brownian motion $B(a_t)$ in distribution.
For $f\in C_b(K/M)$, $f(B(t))-\int_0^t\frac{1}{2}\Delta_{K/M} f(B(s))ds$ is a martingale,
and hence, $f(B(a_t))-\int_0^t\frac{1}{2}\Delta_{K/M}f(B(a_s))da_s$ is a martingale.
On the other hand,  for the conditioned process $z_t$,
$f(z_t)-\int_0^t\frac{1}{2}\Delta_{K/M}f(z_s)da_s$ is a martingale.
The uniqueness of the process with the given representation $(b,A,\Pi)$
implies $z_t=B(a_t)$ in distribution. \ $\Box$

\section{A class of $K$-invariant Markov processes} \label{ex}

In this section, we will study a type $K$-invariant Markov process $x_t$ in $X$ which
is obtained from a $K$-invariant diffusion process interlaced with jumps.

We first consider a $K$-invariant diffusion process $x_t^0$ in $X$. Its generator
takes the following form in local coordinates $x_1,\ldots,x_n$ on $X$:
\begin{equation}
L^0 = \frac{1}{2}\sum_{i,j=1}^n c_{ij}(x)\frac{\partial^2}{\partial x_i\partial x_j} +
\sum_{i=1}^n c_i(x)\frac{\partial}{\partial x_i} \label{L}
\end{equation}
for some $c_{ij},c_i\in C^\infty(X)$ with $c_{ij}$ forming
a nonnegative definite symmetric matrix.

At $x=zy\in X^\circ$, we may assume
$y_1=x_1,\ldots,y_q=x_q$ form local coordinates in $Y^\circ$ around $y$
and $z_1=x_{q+1},\ldots,z_p=x_n$ form local coordinates in $Z=K/M$ around $z$,
where $n=p+q$, $q=\dim(Y)$ and $p=\dim(Z)$. Recall that the radial part $L^{0,Y}$ of
the operator $L^0$ is
given by $(L^{0,Y}f)\circ J=L^0(f\circ J)$ for $f\in C_c^\infty(Y)$.
This implies that $L^{0,Y}=(1/2)\sum_{i,j=1}^q c_{ij}(y)\partial^2/
(\partial y_i\partial y_j)+\sum_{i=1}^q c_i(y)\partial/\partial y_i$, thus, $c_{ij}, c_i
\in C^\infty(Y^\circ)$ for $i,j=1,2,\ldots,q$.
Note that the sum of mixed partials in $L^0$, $\sum_{i=1}^q\sum_{j=1}^p c_{ij}(y,z)
\partial^2/(\partial y_i\partial z_j)$, is independent of the choice
of local coordinates $y_1,\ldots,y_q$ on $Y^\circ$ and $z_1,\ldots,z_p$ on $Z$.

At least locally, there are (smooth) vector fields $\zeta_0,\zeta_1,\ldots,\zeta_q$
on $Y^\circ$ such that $L^{0,Y}=(1/2)\sum_{i=1}^q\zeta_i\zeta_i+\zeta_0$.
We will assume $L^{0,Y}$ is nondegenerate on $Y^\circ$ in the sense that the symmetric
matrix formed by its second order coefficients $c_{ij}(y)$, $i,j=1,2,\ldots,q$,
is positive definite at all $y\in Y^\circ$.
This is independent of choice for local coordinates.
Then $\zeta_1,\ldots,\zeta_q$
are linearly independent at each point of $Y^\circ$, and hence
the sum of mixed partials in $L^0$, $\sum_{i=1}^q\sum_{j=1}^p c_{ij}(y,z)
\partial^2/(\partial y_i\partial z_j)$, may be written as $\sum_{i=1}^q\eta_i\zeta_i$
for uniquely determined vector fields $\eta_1,\ldots,\eta_q$ on $Z$ that may depend
on $y\in Y$. The $K$-invariance of $L^0$ and the linear independence of $\zeta_1,
\ldots,\zeta_q$ imply that $\eta_i$ are $K$-invariant vector fields on $Z$. We have
$L^0 = (1/2)\sum_{i=1}^q(\zeta_i+\eta_i)^2 + \zeta_0 + L'$
for some $K$-invariant second order differential operator $L'$ on $Z$ that may depend
on $y\in Y$.

We may choose vector fields $\theta_1,\ldots,\theta_q,\theta_{q+1},\ldots,\theta_{q+p},
\theta_0$ on $Y^\circ\times Z$ such that $L^0=(1/2)\sum_{i=1}^{q+p}\theta_i\theta_i+
\theta_0$. Let $(\theta_i)_Y$ and $(\theta_i)_Z$ be respectively the components of
$\theta_i$ tangent to $Y^\circ$ and $Z$. After an orthogonal transformation of
$\theta_1,\ldots,\theta_{q+p}$, which may depends on $(y,z)$, we may assume
that $(\theta_i)_Y=0$ for $i=q+1,\ldots,q+p$. Then $(1/2)\sum_{i=1}^q
(\theta_i)_Y(\theta_j)_Y+(\theta_0)_Y=(1/2)\sum_{i=1}^q\zeta_i\zeta_j+\zeta_0$
for $i,j=1,2,\ldots,q$. It follows that after an orthogonal transformation
of $\theta_1,\ldots,\theta_q$, which may depend on $y$,
we may assume that $(\theta_i)_Y=\zeta_i$ for $1\leq i\leq q$.
Then by the discussion in the previous paragraph, $(\theta_i)_Z=\eta_i$ for $1\leq i\leq q$.
It follows that the second order terms in $L'$ are the same as those in
$(1/2)\sum_{i=q+1}^{q+p}\theta_i\theta_i$ and hence the coefficient matrix of the second
order terms in $L'$ is nonnegative definite.

Let $\myk$ and $\mym$ be respectively the Lie algebras of $K$ and $M$, and
let $\myp$ be an $\Ad(M)$-invariant subspace of $\myk$ complementary to $\mym$. Choose
a basis $\xi_1,\ldots,\xi_p$ of $\myp$. As discussed in section~\ref{nonh},
$\Ad(M)$-invariant vectors in $\myk$ may be identified with $K$-invariant vector fields
on $Z=K/M$, and as a $K$-invariant second order differential operator,
$L'=(1/2)\sum_{i,j=1}^p a_{ij}\xi_i\xi_j+\xi_0$ for an $\Ad(M)$-invariant
matrix $a_{ij}$ (symmetric and nonnegative definite) and
an $\Ad(M)$-invariant $\xi_0\in\myp$, both may depend on $y\in Y$.
The generator $L^0$ may now be written as
\begin{equation}
L^0 = \frac{1}{2}\sum_{i=1}^q(\zeta_i+\eta_i)^2 + \frac{1}{2}\sum_{i,j=1}^p
a_{ij}\xi_i\xi_j + \zeta_0 + \xi_0 = L^{0,Y} + L^{YZ} + L_1^Z + L_2^Z \label{L2}
\end{equation}
where $L^{0,Y}=(1/2)\sum_{i=1}^q\zeta_i^2+\zeta_0$ is the radial part of $L^0$,
$\eta_1,\ldots,\eta_q,\xi_0$ are $K$-invariant vector fields on $Z$ which may depend
on $y\in Y$, $L^{YZ}=\sum_{i=1}^q\eta_i\zeta_i$ is the part of $L^0$ containing mixed partials,
$L_1^Z=(1/2)\sum_{i=1}^q\eta_i^2+\xi_0$ and $L_2^Z=(1/2)\sum_{i,j=1}^p a_{ij}(y)\xi_i\xi_j$.
The components $L^{YZ}$, $L_1^Z$ and $L_2^Z$ of operator $L^0$ are
uniquely determined. To show this, note that a different choice for $\zeta_1,\ldots,\zeta_q$
leads to an orthogonal transform of these vector fields and the transposed
transform of $\eta_1,\ldots,\eta_q$. This implies that $L^{YZ}$ and $L_1^Z$, and hence $L_2^Z$,
are not changed. We may write $L_2^Z(y)$ for $L_2^Z$ to
indicate its dependence on $y\in Y$ and similarly for $L^{YZ}$ and $L_1^Z$.

Let $\sigma_{ij}$ be the square root matrix of $a_{ij}$,
that is, $a_{ij}=\sum_{k=1}^p\sigma_{ki}\sigma_{kj}$. Then $L_2^Z=(1/2)\sum_{i=1}^p
(\sum_{j=1}^p\sigma_{ij}\xi_j)^2$.
The $K$-invariant vector fields $\eta_1,\ldots,\eta_q,\xi_0$ on $Z=K/M$ are the
projections of left invariant vector fields $\hat{\eta}_1,\ldots,\hat{\eta}_q,\hat{\xi}_0$
on $K$ which are $\Ad(M)$-invariant at $e$. Replacing $\eta_i$ and $\xi_0$ by
$\hat{\eta}_i$ and $\hat{\xi}_0$, then $L^0$ given by (\ref{L2}) becomes
\[\hat{L} = (1/2)\sum_{i=1}^q(\zeta_i+\hat{\eta}_i)^2 + (1/2)\sum_{i=1}^p
(\sum_{j=1}^p\sigma_{ij}\xi_j)^2 + \zeta_0 + \hat{\xi}_0.\]
This is a differential operator on $Y^\circ\times K$ and is the generator of
the diffusion process $(y_t^0,k_t)$ that solves the
following sde (stochastic differential equation) of Stratonovich form on $Y^\circ\times K$.
\begin{equation}
\left\{\begin{array}{lll} dy_t^0 &=& \sum_{i=1}^q\zeta_i(y_t^0)\circ dB_t^i +
\zeta_0(y_t^0)dt \\
 dk_t &=& \sum_{i=1}^q\hat{\eta}_i(y_t^0,k_t)\circ dB_t^i +
\hat{\xi}_0(y_t^0,k_t)dt +
 \sum_{i=1}^p[\sum_{j=1}^p\sigma_{ij}(y_t^0)\xi_j(k_t)]\circ dB_t^{q+i}
\end{array} \right. \label{sdeyk}
\end{equation}
with $y_0^0=y$ and $k_0=e$, where $B_t=(B_t^1,\ldots,B_t^n)$ is an $n$-dim standard
Brownian motion. Because, for $f\in C_c^\infty(Y^\circ\times Z)$, $(Lf)\circ
(\id_Y\times\pi)=\hat{L}[f\circ(\id_Y\times\pi)]$, where $\pi$: $K\to Z=K/M$ is
the natural projection, it follows that $x_t^0=k_ty_t^0$ with radial part $y_t^0$
and angular part $z_t^0=k_to$. The radial part $y_t^0$ is assumed to be nondegenerate
in the sense that its generator $L^{0,Y}$ is nondegenerate on $Y^\circ$.

Let $k_t=u_ta_t$. Suppose $u_t$ and $a_t$, with $a_0=u_0=e$, solve
the following sde on $K$.
\begin{equation}
\left\{\begin{array}{lll} da_t &=& \sum_{i=1}^q\hat{\eta}_i(y_t^0,a_t)\circ dB_t^i +
\hat{\xi}_0(y_t^0,a_t)dt \\
 du_t &=& \sum_{i=1}^p\{\sum_{j=1}^p\sigma_{ij}(y_t^0)[\Ad(a_t)\xi_j](u_t)
\}\circ dB_t^{q+i}, \end{array} \right. \label{sdeua}
\end{equation}
Then, writing $k\xi$ or $\xi k$ for left or right translations of a left invariant
vector field $\xi$ on $K$ by $k\in K$, and noting $k\xi(k')=kk'\xi(e)=\xi(kk')$ and
$[\Ad(k)\xi](k')=k'k\xi(e)k^{-1}$, we have
\begin{eqnarray*}
&& dk_t = (\circ du_t)a_t + u_t\circ da_t \\
&=& \sum_{i=1}^p\{\sum_{j=1}^p\sigma_{ij}(y_t^0)[\Ad(a_t)\xi_j](u_t)\}a_t\circ dB_t^{q+i} +
\sum_{i=1}^q u_t\hat{\eta}_i(y_t^0,a_t)\circ dB_t^i + u_t\hat{\xi}_0(y_t^0,a_t)dt \\
&=& \sum_{i=1}^p\{\sum_{j=1}^p\sigma_{ij}(y_t^0)\xi_j(u_ta_t)\}\circ dB_t^{q+i} +
\sum_{i=1}^q \hat{\eta}_i(y_t^0,u_ta_t)\circ dB_t^i + \hat{\xi}_0(y_t^0,u_ta_t)dt.
\end{eqnarray*}
Therefore, the second equation in (\ref{sdeyk}) is equivalent to (\ref{sdeua}).

Because $\hat{\eta}_i$ and $\hat{\xi}_0$ are left invariant vector fields on $K$
that are also right $M$-invariant, $a_t$ in the first equation of (\ref{sdeua})
may be replaced by either $a_tm$ or $ma_t$ for any $m\in M$. The uniqueness
of solution implies that $a_tm=ma_t$ for all $m\in M$.
Therefore, $b_t=a_to$ is an $M$-invariant continuous process
in $Z$ with $b_0=o$. Let $v_t=u_to$. Then $z_t^0=u_ta_to=v_tb_t$.

By the first equations in (\ref{sdeyk}) and (\ref{sdeua}), and the nondegeneracy of $L^{0,Y}$,
the $\sigma$-algebras generated by the radial process $y_t^0$ and
by the $q$-dim Brownian motion $(B_t^1,\ldots,B_t^q)$ are the same and contain
the $\sigma$-algebra generated by $a_t$.
Thus, given radial process $y_t^0$, $a_t$ and hence $b_t$ become non-random.
By the second equation in (\ref{sdeua}),
the It\^{o} formula, and the independence between $(B_t^1,\ldots,B_t^q)$ and
$(B_t^{q+1},\ldots,B_t^n)$, for $f\in C_c^\infty(K)$,
\begin{equation}
f(u_t) - \int_0^t(1/2)\sum_{i,j=1}^pa_{ij}(y_s^0)[\Ad(a_s)\xi_i][\Ad(a_s)\xi_j]f(u_s)ds
\label{fut}
\end{equation}
is a martingale given process $y_t^0$. Then, replacing $u_t$ and $a_t$ by $v_t$ and $b_t$,
the above is still a martingale for $f\in C_c^\infty(Z)$ given process $y_t^0$.
Recall $\zeta$ is the exit time of the radial process from $Y^\circ$. The above shows that
almost surely on $[\zeta>T]$, given $y_t^0$ for $0\leq t\leq T$, the conditioned
angular process $z_t^0$, as a continuous nonhomogeneous \levy process
in $Z=K/M$, has drift $b_t$ and covariance function $A_{ij}(t)=\int_0^t a_{ij}(y_s^0)ds$
under basis $\xi_1,\ldots,\xi_p$ on $\myp$. Its covariance operator
is thus $\int_0^t L_2^Z(y_s^0)ds=\int_0^tds\frac{1}{2}\sum_{i,j}a_{ij}(y_s^0)\xi_i\xi_j$.

Let $G$ be a topological group acting on $X$ continuously and containing $K$ as
a topological subgroup. A measure $\mu$ on $G$
is said to be $K$-conjugate invariant if $c_k\mu=\mu$ for $k\in K$,
where $c_k$ is the conjugation map: $G\ni g\to kgk^{-1}\in G$.

Now let $x_t$ be a process in $X$ obtained from a $K$-invariant diffusion
process $x_t^0$ in $X$ interlaced with jumps determined by a Poisson random measure $N$
on $\bR_+\times G$ of intensity measure $dt\eta(d\sigma)$, where $\eta$ is
a finite $K$-conjugate invariant measure on $G$. The process $x_t$ may be constructed
as follows.
Let $\gamma_n$ be positive random variables of a common exponential distribution of
mean $1/\eta(G)$ and let $\sigma_n$ be random variables in $G$
of common distribution $\eta/\eta(G)$, all independent of each other and of process $x_t^0$,
and let $T_n=\gamma_1+\cdots+\gamma_n$ with $T_0=0$.
Then the process $x_t$ is obtained by setting $x_t=x_t^0$ for $t<T_1$ and
inductively setting $x_t=x_t^n$ for $T_n\leq t<T_{n+1}$, where $x_t^n$
is $x_t^0$ determined by the initial condition $x_{T_n}^n=\sigma_n(x_{T_n-})$.

It is easy to see that $x_t$ defined above is a Markov process in $X$ and
its transition semigroup $P_t$ is the solution of the following integral equation:
for $f\in C_0(X)$,
\begin{equation}
\forall x\in X,
\ \ \ \ P_tf(x) = e^{-\lambda t}P_t^0f(x) + \int_0^t e^{-\lambda u}du
\int_G\int_XP_u^0(x,dx_1)P_{t-u}f(\sigma x_1)\eta(d\sigma) \label{inteq}
\end{equation}
with $P_0f=f$, where $\lambda=\eta(G)$ and $P_t^0$ is the transition semigroup of $x_t^0$
(a Feller semigroup). The existence of the solution $P_tf$ and the Feller property
of $P_t$ can be established by a routine argument of successive approximation,
and the uniqueness of solution can be shown by a simple estimate on the difference
between two possible solutions. By the uniqueness of $P_t$, the $K$-invariance
of $P_t^0$ and the $K$-conjugate invariance of $\eta$, it is easy to show that $P_t$
is $K$-invariant. Thus, $x_t$ is a $K$-invariant Feller process in $X$.
Moreover, by differentiating (\ref{inteq}) at $t=0$, we obtain its generator $L$
on $f\in C_c^\infty(X)$:
\begin{equation}
Lf(x) = L^0f(x) + \int_G[f(\sigma x)-f(x)]\eta(d\sigma), \label{L3}
\end{equation}
where $L^0$ is the generator of the $K$-invariant diffusion process $x_t^0$.

We need some notation in order to describe the \levy measure function of
the conditioned angular process of $x_t$.
Recall $J$ and $J_2$ are the projection maps $X\to Y$ and $X^\circ\to Z=K/M$.
For a finite measure $\mu$ on $X$ and $y\in Y$, let $\mu(\cdot\mid y)$ be
the conditional distribution of $\mu$ given the orbit $Ky$, which is the
regular conditional distribution of a random variable $x$ with distribution $\mu/\mu(X)$
given $Jx=y$ and is a probability measure supported by $Ky$. 
Note that if $\mu(Ky)>0$, then $\mu(\cdot\mid y)=\mu(\cdot\,;Ky)/\mu(Ky)$,
where $\mu(\cdot\,;Ky)$ is the restriction of $\mu$ to $Ky$, but even if $\mu(Ky)=0$,
$\mu(\cdot\mid y)$ is still defined for $J\mu$-almost all $y$.
If $y\in Y^\circ$, then $\mu(\cdot\mid y)$ is supported by $X^\circ$ and hence $J_2$ may
be applied to obtain a probability measure $J_2\mu(\cdot\mid y)$ on $Z$.
For $x\in X$, let $i_x$ be the evaluation map: $G\ni g\mapsto gx\in X$
and for any measure $\eta$ on $G$, let $\eta(x)$ be the measure $i_x\eta$ on $X$,
that is, $\eta(x)(f)=\eta(f\circ i_x)$ for $f\in C_b(X)$.
Note that if $\eta$ is $K$-conjugate invariant, then $k\eta(x)=\eta(kx)$ for $k\in K$,
and hence, $\eta(y)$ is $M$-invariant for $y\in Y^\circ$.
It then follows that $\eta(y_1)(\cdot\mid y_2)$ is $M$-invariant
and hence $J_2\eta(y_1)(\cdot\mid y_2)$ is an $M$-invariant probability measure on $Z$.

The points in the Poisson random measure $N$ are the pairs
$(\tau,\sigma)$, where $\tau>0$ is a jump time and $\sigma\in G$ is the associated
jump size of process $x_t$ in the sense that $x_\tau=\sigma x_{\tau-}$.
Under the decomposition $x=zy$ for $x\in X^\circ$ with $y=J(x)$ and $z=J_2(x)$,
\begin{equation}
z_\tau y_\tau = \sigma z_{\tau-}y_{\tau-} = S(z_{\tau-})\sigma'y_{\tau-} = S(z_{\tau-})
[J_2(\sigma' y_{\tau-})][J(\sigma' y_{\tau-})], \label{zytau}
\end{equation}
where $\sigma'=S(z_{\tau-})^{-1}\sigma S(z_{\tau-})$ and $S$ is a section map
on $Z=K/M$. Because $\sigma$ is independent of $x_{\tau-}$ and
its distribution $\eta/\eta(G)$ is $K$-conjugate invariant, $\sigma'$ has
the same distribution as $\sigma$ and is independent of $x_{\tau-}$,
and hence $(\sigma' y_{\tau-})$ has distribution
$\eta(y_{\tau-})(\cdot)/\eta(y_{\tau-})(X)$. By (\ref{zytau}) and the uniqueness of
the decomposition $x=zy$, $J(\sigma' y_{\tau-})=y_\tau$ and $S(z_{\tau-})
[J_2(\sigma' y_{\tau-})]=z_\tau$. Thus, the jump of process $z_t$
at time $\tau$ is $z_{\tau-}^{-1}z_\tau=S(z_{\tau-})^{-1}z_\tau=J_2(\sigma' y_{\tau-})$.
Here the product on $Z=K/M$ should be understood in the sense described in section~\ref{nonh},

Note that almost surely, $x_t$ and hence $y_t$ have finitely many jumps
for $0\leq t\leq T$. Given $y_t$ in $Y^\circ$ for $0\leq t\leq T$ with jumps at
fixed times $t_1<t_2<\cdots<t_k$ means three things:
First, the Poisson random measure $N$ is conditioned to
have points $(t_i,\sigma_i)$ for $1\leq i\leq k$ with $\sigma_i$ satisfying
$J(S(z_{t-})^{-1}\sigma_i S(z_{t-})y_{t-})=y_t$ for $t=t_i$; secondly, at all other     
points $(\tau,\sigma)$ of $N$, $J(S(z_{\tau-})^{-1}\sigma S(z_{\tau-})y_\tau)=y_\tau$; 
and finally, $y_t^0$ is given on each time interval $[t_{i-1},\,t_i)$.
Note that when $N$ is conditioned to have points $(t_i,\sigma_i)$ for
$1\leq i\leq k$, its remaining points form a Poisson random measure $N'$ of same distribution
as $N$. 
On the time interval $[t_{i-1},\,t_i)$, the conditioned angular process $z_t$
may be obtained from the diffusion process $z_t^0$ by adding jumps from $N'$.
It follows that given $y_t$ in $Y^\circ$ for $0\leq t\leq T$, the conditioned
process $z_t$, as a nonhomogeneous \levy process in $Z=K/M$, has possible fixed
jumps $z_{\tau-}^{-1}z_\tau$, at $\tau=t_1,t_2,\ldots,t_k$, of distribution
$J_2\eta(y_{\tau-})(\cdot\mid y_\tau)$ (this measure may be equal to $\delta_o$ and then
there is no fixed jump at $\tau$), its remaining jumps form a Poisson random measure
$N^Z$ on $[0,\,T]\times Z$ with intensity measure
$\mu(dt,dz)=dt[J_2\eta(y_t)(\cdot\,;Ky_t)](dz)$. The \levy measure function of $z_t$ is
$\Pi(t,\cdot)=\int_0^t\mu(ds,\cdot)$. Its covariance operator
is the same as that of $z_t^0$ obtained earlier (but with $y_t^0$ replaced by $y_t$).
To summarize, we have the following result.

\begin{theo} \label{th5}
Let $x_t$ be a $K$-invariant Markov process in $X$ obtained from
a $K$-invariant diffusion process $x_t^0$ interlaced with jumps given by
a Poisson random measure $N$ on $\bR_+\times G$ with intensity measure $dt\eta(d\sigma)$
as described above,
where $G$ is a topological group acting continuously on $X$ and containing $K$
as a topological subgroup, and $\eta$ is a finite $K$-conjugate invariant
measure on $G$. Then $x_t$ is
a Feller process in $X$ with generator given by (\ref{L3}). 

Assume $x_t^0$ has a nondegenrate radial part and $x_0\in X^\circ$. Let $x_t=z_ty_t$ be
the decomposition into radial process $y_t$ in $Y^\circ$ and angular process
$z_t$ in $K/M$ for $t<\zeta$, where $\zeta$ is the exit time of
$x_t$ from $X^\circ$ or $y_t$ from $Y^\circ$, and let $L_2^Z(y)$ be the part of
generator $L^0$ of $x_t^0$ given in (\ref{L2}). Fix $T>0$. Then almost surely on
$[\zeta>T]$, given $y_t$ for $0\leq t\leq T$ with jumps at $t_1<t_2<\cdots<t_k$,
the conditioned process $z_t$ is a nonhomogeneous \levy process in $K/M$
with possible fixed jumps $z_{\tau-}^{-1}z_\tau$ at $\tau=t_1,t_2,\ldots,t_k$
of distribution $J_2\eta(y_{\tau-})(\cdot\mid y_\tau)$, covariance operator
$V(t)=\int_0^t L_2^Z(y_s)ds$ and \levy measure function
$\Pi(t,dz)=\int_0^t\,ds\,[J_2\eta(y_s)(\cdot\,;Ky_s)](dz)$.
\end{theo}

\end{document}